\documentclass{amsart}
\usepackage{amsfonts,amscd, amssymb}

\newtheorem{theorem}{Theorem}[section]
\newtheorem{lemma}[theorem]{Lemma}

\newtheorem{proposition}[theorem]{Proposition}
\theoremstyle{remark}

\theoremstyle{definition}
\newtheorem{definition}[theorem]{Definition}

\numberwithin{equation}{section} \makeatother

\begin{document}

\title[Closed projections for operator algebras]{Closed projections and peak interpolation for operator algebras} 

\author{Damon M. Hay}
\address{Department of Mathematics, University of Houston, Houston, TX
  77204-3008}
\email[Damon Hay]{dhay@math.uh.edu}  

\begin{abstract}The closed one-sided ideals of a $C^*$-algebra are
exactly the closed subspaces supported by the orthogonal complement
of a closed  projection.  Let $A$ be a (not necessarily selfadjoint)
subalgebra of a unital $C^*$-algebra $B$ which contains the unit of
$B$.  Here we characterize the
right ideals of $A$ with left contractive approximate identity as
those subspaces of $A$ supported by the 
orthogonal complement of a closed
projection in $B^{**}$ which also lies in $A^{\perp \perp}$.  Although
this seems quite natural, the proof requires a set of new techniques which
may may be viewed as a noncommutative version of the subject of peak
interpolation from the 
theory of function spaces.  Thus, the right ideals with left approximate identity are closely related to a type of peaking phenomena in the algebra.  In this
direction, we introduce a class of closed projections 
which generalizes the notion of a peak set in the theory of uniform
algebras to the world of operator algebras and operator spaces.
\end{abstract}  

\maketitle

\section{Introduction}Let $K$ be a compact Hausdorff space and let
$C(K)$ denote the $C^*$-algebra of all complex-valued continuous
functions on $K$.  It is well known that closed ideals in $C(K)$
consist of all functions which vanish on a fixed closed subset of $K$.
If instead, $A$ is a uniform algebra contained in $C(K)$, then by a theorem of 
Hirsberg \cite{Hirsberg}, a closed subspace $J$ of $A$ is a closed 
ideal with contractive approximate identity if and only if it consists of all 
functions which vanish on a `p-set' 
for $A$.  Recall, a subset $E$ of $K$ is said to be a \emph{peak set} for $A$ 
if there 
exists a function $f$ in $A$ such that $f(x) =1$ for all $x \in E$ and 
$|f(x)| < 1$ for all $x \in E^c$.  A subset $E$ of $K$ is said to be a 
\emph{p-set} for $A$ if it is the intersection of a family of peak sets.
The p-sets for a uniform algebra $A$ were characterized by Glicksberg
\cite{G} as those 
closed subsets $E$ such that $\mu \in A^\perp$ implies $\mu_E \in A^\perp.$  See \cite{Gam}, \cite{G}, or \cite{Jar} for more information on peak 
sets.

For general $C^*$-algebras, the closed right ideals of a
$C^*$-algebra $A$ consist of the elements $a$ in $A$ for which $qa = 0$ for a
\emph{closed} projection $q$ in the second dual of $A$.  In other
words, a subspace $J$ is a right ideal of $A$ if and only if $$J =
(1-q)A^{**} \cap 
A,$$ for a closed projection $q$ in $A^{**}$.  Of course we are viewing
$A$ as being canonically embedded in its second dual, which is a
$W^*$-algebra.  In fact, $1-q$ will be a weak*-limit point for any
left contractive approximate identity of $J$.  Indeed all closed projections
arise in this manner.

Turning to the nonselfadjoint case, let $A$ be a subalgebra of a
unital $C^*$-algebra $B$, such that $A$ contains the identity of $B$.
We characterize the right ideals of $A$ with left contractive
approximate identity as those subspaces $J$ of the form $J = (1-q)A^{**}
\cap A$, for a closed (with respect to $B^{**}$) projection $q$ in
$A^{**}$.  However natural 
this may appear, the tools available in the selfadjoint theory are not
applicable here.  Thus a portion of this paper develops
some technical 
tools from which this characterization follows.  Incidentally, these
generalize some peak interpolation results in the theory of uniform
algebras.  
The above mentioned characterization is a refinement of the characterization in
\cite{BEZ}, which is in terms of right $M$-ideals.  In particular, it
appears to open up
a new area in the theory of nonselfadjoint operator algebras,
allowing for the generalization of certain important parts of the theory of
$C^*$-algebras.  This will be explored more fully in the sequel \cite{BHN} 
where, for example, we apply the main result of this paper to develop a theory 
of hereditary subalgebras of not necessarily selfadjoint operator algebras. 
As is the case in the selfajoint theory, we demonstrate that these hereditary 
subalgebras are connected to the facial structure of the state space.  
Additionally, we also give a solution there to a more than ten year old problem
 in the theory of operator modules.

In our noncommutative setting, the peak and $p$-sets described above
are replaced with a certain class of projections in the second dual of $B$,
called the \emph{peak} or \emph{$p$-projections for $A$}.  In the commutative
case this class of projections  
can be identified with the characteristic functions of peak or
$p$-sets for $A$.  
When $A = B$, the $p$-projections are exactly
the closed projections in $A^{**}$.  The theory of these projections brings 
another tool from the classical theory to the world of operator spaces.

The paper is organized as follows.  In Section \ref{prelim} we
introduce the notation and discuss some background and preliminary
results.  In particular, we discuss the 
noncommutative topology of open and closed projections.  Section
\ref{interpolation} generalizes some interpolation results from the
theory of function spaces to operator spaces, which will be used in
Sections \ref{right} and \ref{peak}.  Section \ref{right} contains the
main theorem and its proof.  Finally, in Section \ref{peak} we look at
closed projections in the weak*-closure of an operator algebra from
the perspective of `peak phenomena.'

\bigskip

{\bf Acknowledgments.} \ We first thank our Ph.D. advisor, David Blecher, for 
the initial inspiration to pursue this project and for his
continuous support and suggestions along the way. We also thank him
for pointing out, and helping to correct, some errors in the proof of Proposition
\ref{gamlemma1} in an earlier version.  Additionally, we extend our thanks to
Professor Charles Akemann for clarifying a number of points regarding
open and closed projections.

\section{Open and closed projections and preliminary results}\label{prelim}
The theory of operator spaces and completely bounded maps has
long been recognized as the appropriate setting for studying many
problems in operator algebras.  Basics on operator spaces may be
found in \cite{BL}, \cite{ER}, \cite{Paulsen}, and \cite{Pisier}.  We
will make use of the following lemma which gives a criteria for when a
completely bounded map is a complete isometry.

\begin{lemma}\label{ci} Let $X$ be an operator space and $Y$ a (not necessarily
  closed) subspace of another operator space.  Suppose $T:X \to
  Y$ is a one-to-one and surjective completely bounded map such that
  $T^*$ is a complete isometry.  Then $T$ is a complete isometry. 
\end{lemma}

\begin{proof}Let $Z$ be the closure of $Y$ and define $R$ to be the
  same as $T$ except with range $Z$.  Since $Z^* = Y^*$, we have that
  $R^* : Z^* \to X^*$ is simply $T^*$ which is one-to-one and has
  closed range.  By VI.6.3 in \cite{Dun} this implies that $R$ is
  onto.  Since $R$ 
  is onto, $Y = Z$ and $R = T$, and by the open mapping theorem $T$ is
  bicontinuous.  Hence, given $\varphi \in  X^*$, then $\psi \equiv \varphi
  \circ T^{-1} \in Y^*$ satisfies $T^*\psi = \varphi$, showing that $T^*$ is
  surjective and thus, since it is also completely isometric, $T^{**}$
  is a complete isometry. 
  Viewing $X$ and $Y$ as being canonically embedded in their second
  duals, $T$ is just the restriction of $T^{**}$ to $X$.  Hence, $T$
  is also a complete isometry.
\end{proof}

Throughout this paper, $B$ will denote a unital $C^*$-algebra and $A$
will denote either a unital subspace or a unital subalgebra of $B$, for
which by \emph{unital} we mean $1_B \in A$.  We
view $A$ and $B$ as being canonically embedded into the second dual
$B^{**}$ of $B$ via the canonical isometry.  The second dual of $B$,
$B^{**}$, is a $W^*$-algebra.  By the \emph{state space of $B$}, which
we denote $S(B)$, we
mean the set of positive functionals on $B$ which have norm one.
Each functional $\varphi$ of $B$ extends uniquely to a weak*-continuous, or
\emph{normal}, functional on $B^{**}$, which we again denote by $\varphi$.
Also, we denote the unit in $B$ by $1_B$, or more often, simply by $1$.

By a \emph{projection} in $B$ or $B^{**}$ we mean an orthogonal projection.  
The \emph{meet} of any two projections $p$ and $q$ can be given
abstractly as $$p \wedge q = \lim_{n \rightarrow \infty}(pq)^n,$$
where this limit is taken in the weak* topology.  Similarly, the join
is given by 
$$p \vee q = \lim_{n \rightarrow \infty}1-(1-p-q+pq)^n.$$
Let $M$ be a (not necessarily selfadjoint) weak*-closed unital
subalgebra of a 
$W^*$-algebra and suppose that $p$ and $q$ are projections in
$M$.  Then by the formula for $p \wedge q$ above, $p \wedge q$ is also in $M$.  By
induction, this extends to the meet of any finite 
collection of projections $M$.  More generally, if $\{p_{\alpha}\}$ is
any collection of projections, then $\wedge p_{\alpha}$ is the weak*-limit of
the net of meets of finite subcollections of $\{p_{\alpha}\}$, each
of which is in $M$.  Thus $\wedge p_{\alpha}$ is also in $M$.
Similarly, a join of projections in $M$ is also in $M$.  Now let $\pi
: M \rightarrow B(H)$ be a weak*-continuous homomorphism of $M$ into
the bounded linear operators on a Hilbert space $H$. If $p$ and $q$
are projections in $M$, then $\pi((pq)^n) = (\pi(p)\pi(q))^n$ for each
$n$, and so we have $\pi (p \wedge q) = \pi(p) \wedge
\pi(q)$.  This clearly generalizes to meets of finitely many
projections.  By approximating by such finite meets, this in turn
generalizes to arbitrary meets of projections.  Similar statements
apply to joins of projections.


A projection $p \in B^{**}$ is said to be \emph{open} if it is the
weak*-limit of an increasing net $(b_t)$ of elements in $B$ with $0
\le b_t \le 1$.  A projection $q \in B^{**}$ is said to be
\emph{closed} if $1-q$ is open. It is clear that a closed projection
is the weak*-limit of a decreasing net of positive elements in $B$.
It is well known that a projection $p$ in $B^{**}$ is open if and
only if it
is the support of a left (respectively, right) ideal in $B$.  That is,
there exists 
a left (respectively, right) ideal $J$ in $B$ such that $J = B^{**}p
\cap B$ (respectively, $J = pB^{**}\cap B$).  In this case, 
the weak* closure of $J$ in $B^{**}$ is $B^{**}p$ (respectively,
$pB^{**}$).  Moreover, $p$ is a 
weak*-limit point of any increasing right contractive approximate
identity for $J$.  Actually, if $p \in B^{**}$ is a projection which is
a weak* limit of a net $(e_t)$ in $B$ such
that $e_tp = e_t$, then $p$ is open. To see this, we let $J$ be the
set of all $b \in B$ such that $bp = b$.  Then $J$ contains $(e_t)$
and so $p$ is in $J^{\perp \perp}$, which is a weak*-closed left ideal
of $B^{**}$. Thus $B^{**}p \subset J^{\perp \perp}$, but also
$J^{\perp \perp} \subset B^{**}p$, so that $J^{\perp \perp} =
B^{**}p$.   However, $J = B^{**}p \cap B$, so that  $p$ is the support
of a closed left ideal, making it 
an open projection.  A similar argument using right ideals holds if
$pe_t = e_t$ instead. 
In the case that $B$ is commutative, open and closed projections correspond to
characteristic functions of open and closed sets, respectively.
It is this collection of open and closed projections which will act as a
kind of substitute for topological arguments in the noncommutative
situation.  We now list, most without proof, some basic facts
regarding these open and closed projections.  Many of these facts can
be found in Akemann's papers \cite{Ak1} and \cite{Ak2}, and some may
also be found in \cite{Ped} and \cite{GK}. 

The join of any collection of open projections is again an open
projection.  Hence, the meet of any collection of closed projections
is again a closed projection.  However, in contrast to the commutative
situation, joins of closed projections are not necessarily closed (see
\cite{Ak1}).   For a general $C^*$-algebra $B$,
the presence of a unit guarantees a kind of noncommutative
compactness.  That is, if $q$ is a closed projection,  
given any collection of open projections $\{p_\alpha\}$ such that
$q \le \bigvee_\alpha p_\alpha$, then there exists a finite subcollection
$\{p_{\alpha_1}, \dots ,p_{\alpha_k}\}$ such that
$q \le \bigvee_{i=1}^k p_{\alpha_i}$ (see Proposition II.10 in
\cite{Ak1}).  We will refer to this as the `compactness property.'  
A type of regularity also holds with respect to open and closed
projections.  Namely, any closed 
projection is the meet of all open projections dominating it.  The
following was communicated to us by Akemann.

\begin{proposition}\label{regular}Let $q$ be a closed projection in $B^{**}$.  Then
$$q = \bigwedge\{u\; | \; u \ge q \; \text{and open} \}.$$\end{proposition}

\begin{proof}Assume that $B$ and $B^{**}$ are represented in the universal
    representation of $B$.  
    Since $q$ is closed we may find an increasing net $(a_t)$ in
    $B_{sa}$ such that $1-a_t \searrow q$ weak* and $(1-a_t)q = q$.
    Using the Borel functional calculus, for each $t$ let $$r_t =
    \chi_{(\frac{1}{2} , \infty)}(1-a_t),$$ where $\chi_{(\frac{1}{2} ,  \infty)}$ is the
    characteristic function of the open interval $(\frac{1}{2} ,
    \infty)$.  Then necessarily $r_t$ is an open projection
    such that
    $$2(1-a_t) \ge r_t.$$  Furthermore,
    we claim that each $r_t$ dominates $q$.  To prove this claim,
    fix $t$ and let $\{f_n\}$ be an increasing
    sequence of positive continuous functions on the spectrum of $1-a_t$ which
    converges point-wise to $\chi_{(\frac{1}{2} , \infty)}$ and is such
    that $f_n(1) = 1$ for each $n$.  Now fix $n$ and suppose $\xi
    \in \text{Ran} \; q$ is of norm one.  Then $(1-a_t)\xi =
    (1-a_t)q\xi = q\xi = \xi$.  So for any 
    polynomial $R$, we have $R(1-a_t)\xi = R(1)\xi$.  Let $R_k$ be a
    sequence of polynomials converging uniformly to $f_n$.  Then
    $$\langle f_n(1-a_t)\xi,\xi \rangle = \lim_k \langle R_k(1)\xi,\xi\rangle =
    \lim_kR_k(1) = f_n(1) = 1.$$  By the converse to the
    Cauchy-Schwarz inequality, we have $f_n(1-a_t)\xi = \xi$ for all
    $\xi \in \text{Ran} \; q$, which is to
    say that $f_n(1-a_t)q = q$.  Hence
    $\chi_{(\frac{1}{2},\infty)}(1-a_t)q = q$, and so $r_t \ge q$.
    Let $r_0 =
    \bigwedge_{t}r_t$ and suppose that $r_0$ does
    not equal $q$.  Then there exists a state $\varphi \in S(B)$ such
    that $\varphi(r_0 - q) = 1$.  This forces $\varphi(r_0) = 1$ and
    $\varphi(q) = 0$ since $r_0 - q \ge 0$.  Thus $$\varphi \left(
      2(1-a_t) \right) \rightarrow 0.$$  However, 
    since $\varphi(r_0) = 1$ and $r_t \ge r_0$, it must be
    that $\varphi(r_t) = 1$.  Applying $\varphi$ to the
    inequality $$r_t \le 2(1-a_t)$$ and taking the weak* limit, we get
    $1 \le 0$.  Hence, $r_0 = q$ which proves the result.
\end{proof}

Finally, one of the most important results in basic topology is
Urysohn's lemma.  Akemann has extended this result to closed
projections:

\begin{theorem}[\cite{Ak2}]\label{ncu}Let $p$ and $q$ be closed projections in
  $B^{**}$ for a 
  $C^*$-algebra $B$ such that $pq = 0$.  Then there exists an element
  $a$ in $B$, $0 \le a \le 1$, such that $ap = p$ and $aq = 0$.
\end{theorem}

We will often be working with closed projections in $B^{**}$ which lie
in the weak* closure of $A$ in $B^{**}$.  The following gives some
equivalent conditions for this.

\begin{lemma}\label{weak}
Let $A$ be closed subalgebra of a unital $C^*$-algebra $B \subset B^{**}$ such
that $A$ contains the unit of $B$.  Let $q \in B^{**}$ be a projection.
The following are equivalent:
\begin{enumerate}
 \item $q \in \overline{A}^{w*},$
 \item $q \in A^{\perp \perp},$
 \item $A^{\perp} \subset (qA)_{\perp}.$
\end{enumerate}
\end{lemma}

\begin{proof}
  The equivalence of (1) and (2) is a standard result of functional
  analysis.  Suppose (3) holds.  Then $((qA)_{\perp})^{\perp}
  \subset A^{\perp \perp}$.  However, $((qA)_{\perp})^{\perp} =
  \overline{qA}^{w^*} = q \overline{A}^{w^*}$ which must contain $q$
  since $A$ is unital. Hence, (2) holds.  Now assume (2).  By
  hypothesis, $\psi(q) = 0$ for all $\psi \in 
  A^{\perp}$.  Let $\varphi \in A^{\perp}$.  Then for each $a \in A$,
  $\varphi(\cdot a) \in A^{\perp}.$  Thus $\varphi(qa) = 0$ for all $a
  \in A$.  Hence $\varphi \in (qA)_{\perp}$.
\end{proof}

A class of operators which will play an
important role here are the \emph{completely non-unitary}, or \emph{c.n.u.},
operators on a Hilbert space $H$.  A contraction $T$ is said to be \emph{completely
  non-unitary} if there exists no reducing subspace for $T$ on which $T$ acts
unitarily.  It is well known that if $T$ is completely non-unitary,
then $T^n \rightarrow 0$ in the weak operator topology on $B(H)$ as
$n \rightarrow \infty$.  See \cite{FN} and \cite{Kubrusly} for details.

If $B$ is a unital $C^*$-algebra, we denote the self-adjoint part of
$B$ by $B_{sa}$.  Kadison's `function representation' says that
$B_{sa}$ may be represented as continuous affine functions on $S(B)$
via an order preserving linear isometry which extends
weak*-continuously to $B_{sa}^{**}$, in such a way that
$B_{sa}^{**}$ is represented as bounded affine functions on $S(B)$.
We say that an element $b$ of $B_{sa}^{**}$ is lower semi-continuous
if its image under this representation is a lower semi-continuous
function on $S(B)$ (\cite{Ped}).  

\begin{lemma}\label{vanishing}
 Let $b$ be a positive, lower semi-continuous contraction in $B^{**}$
 for a $C^*$-algebra $B$ and suppose $\varphi_0(b) = 0$ for some
 $\varphi_0 \in S(B)$.  Then there exists a pure state of 
 $B$ which is zero at $b$.
\end{lemma}

\begin{proof}
  Let $K = \{\varphi \in S(B)\; | \; \varphi(b) = 0\}$.  The set $K$
  is nonempty 
  by hypothesis, and since $\phi$ and $b$ are positive, we also have
  that $K  =  \{\varphi \in S(B)\; | \; \varphi(b) \le 0\}$.  Thus $K$
  is the complement of the set $\{\varphi \in S(B)\; | \; \varphi(b) >
  0\}$ which is open in the weak* topology by the semi-continuity
  property of $b$.  Thus $K$ is weak*-closed in $S(B)$ and hence
  weak*-compact.  It is also convex by a straight-forward calculation.
  Thus, $K$ is well supplied with extreme points by the Krein-Milman
  theorem.  Now suppose that $\varphi_1,\varphi_2 \in S(B)$, 
  $\lambda$ is a scalar in $(0,1)$ and that $\lambda\varphi_1 +
  (1-\lambda)\varphi_2 \in K$.  Then $\lambda\varphi_1(b) +
  (1-\lambda)\varphi_2(b) = 0$.  However, by positivity, this forces
  $\varphi_1(b) = \varphi_2(b) = 0$ and so $\varphi_1$ and $\varphi_2$
  are in $K$.  In other words, $K$ is a face of
  $S(B)$ and hence must contain an extreme point of $S(B)$.  However,
  the extreme points of $S(B)$ are the pure states of $B$. 
\end{proof}

\section{Noncommutative peak interpolation}\label{interpolation} 

The following sequence of propositions and lemmas are the keys to the
main result and generalize some classical results from the theory of function 
spaces (see Section II.12 of \cite{Gam}).

\begin{proposition}\label{gamlemma1}Let $X$ be a closed subspace of a
  $C^*$-algebra 
  $B$.  Let $q \in B^{**}$ be a projection such that $\varphi
  \in (qX)_{\perp}$ for all $\varphi \in X^{\perp}$.
  Let $I = \{x \in X : qx = 0\}$.  Then $qX$ is completely isometric
  to $X/I$ via the map $x+I \mapsto qx$. Similarly, if $I$ is defined
  to be $\{x \in X : xq = 0\}$,
  then $Xq$ is completely isometric to $X/I$ via the map $x+I
  \mapsto xq$.
\end{proposition}

\begin{proof}  We will be using standard operator space duality
  theory, as may be found in \cite{BL}, for example.
  First note that $I$ is the kernel of the completely contractive map $x
  \mapsto qx$ on $X$, so
  that this map factors through the quotient $X/I$:
  $$X \overset{S}{\rightarrow} X/I \overset{T}{\rightarrow} qX,$$
  where $S$ is the natural quotient map and $T$ is the induced linear
  isomorphism.  Taking adjoints, we have  $(T \circ S)^* = S^* \circ
  T^*$ and if $\varphi \in (qX)^*$ and $x \in X$, then 
  $$(S^* \circ T^*)(\varphi)(x) = S^*(T^*\varphi)(x)=(T^*\varphi)(Sx)
  = \varphi(TSx) = \varphi(T(x+I))=\varphi(qx),$$
  so that $S^*\circ T^*$ is given by
  $$\varphi \mapsto \varphi(q \cdot ),$$ for each $\varphi \in (qX)^*$.
  Identifying $(qX)^*$ with $(qB)^*/(qX)^{\perp}$ and $X^*$ with
  $B^*/X^\perp$, the map $S^* \circ T^*$  takes an element $\varphi +
  (qX)^\perp$ to the element $\varphi(q \cdot ) + X^\perp$.  To show that
  $T$ is a complete isometry, by Lemma \ref{ci}, it suffices to show
  that $T^*$ is a
  complete isometry, since $T$ is one-to-one, surjective, and
  completely bounded.  Since $S^*$ is completely contractive, if $S^*
  \circ T^*$ is completely isometric, then $$\|\varphi\| = \|(S^*
  \circ T^*)(\varphi)\| \le \|T^*\varphi\| \le \|\varphi\|.$$
  Similar statements also hold for each matrix level.  Hence, if $S^*
  \circ T^*$ is completely isometric, then so is $T^*$.
  Thus, in order to show that $T^*$ is
  completely isometric, it is sufficient to show that $S^* \circ T^*$ is
  a complete isometry.  Note that since $T \circ S$ is completely
  contractive, so is $S^* \circ T^*$.  We let $(e_t)$ be a decreasing
  net in the unit ball of $B$,
  such that $e_t \rightarrow q$ weak* and $qe_t = q$ for all $t$.  
  
  Let $\varphi \in (qB)^*$ and  $\psi \in X^{\perp}$.  Then $\psi \in (qX)_{\perp} \subset
  (qX)^\perp$.  If $J$ is the right ideal in $B$
supported by $1-q$, then for $qb \in {\rm  Ball}(qB)$, we have
$\|qb\| = \|b + J\|$.  Since right ideals are
proximinal in a $C^*$-algebra, it follows that there exists $a \in J$
such that $\|qb\|= \|b+J\| = \|b+a\|$.  Since $q(b+a) = qb$ and
$\|b+a\| \le 1$, by replacing $b$ with $b+a$ it
follows that 
\begin{eqnarray*}\| \varphi + \psi \|_{(qB)^*} & = & \sup \{ |
  \varphi(qb) + \psi (qb)| : b \in {\rm Ball}(B) \}. \end{eqnarray*}

However, for $b \in {\rm Ball}(B)$, we have
\begin{eqnarray*} | \varphi(qb) + \psi (qb)| & = & \lim_t |
  \varphi(qe_tb) + \psi (e_tb)| \\
& \le & \|\varphi (q \cdot) + \psi\|_{B^*} \|e_tb\| \\
& \le & \|\varphi (q \cdot) + \psi\|_{B^*}. \\
\end{eqnarray*}
Hence, $\| \varphi + \psi \|_{(qB)^*} \le \|\varphi (q \cdot) + \psi\|_{B^*},$
and, thus,
$\|\varphi + (qX)^{\perp}\| \le \| \varphi (q \cdot) + \psi\|_{B^*}.$  Now
taking the infimum over all $\psi \in X^\perp$, we get $\|\varphi +
(qX)^{\perp}\| \le \| \varphi (q \cdot) + X^\perp \|$.

The matricial case is almost identical, using operator space duality
principles, and is left to the reader to fill in the details.
The last statement of the
  proposition follows by a completely analogous proof.
\end{proof}

Since $qX \subset qB$ and $qB$ can be identified with a quotient
$B/J$, where $J$ is the right ideal of $B$ corresponding to $q$, then
the result above shows that the set $\{x + J| x \in X\}$ is closed in $B/J$.

\begin{proposition}\label{epsilon}
  Let $X$ be a closed subspace of a $C^*$-algebra
  $B$.  Let $q \in B^{**}$ be a projection such that $\varphi
  \in (qX)_{\perp}$ whenever $\varphi \in X^{\perp}$.
  Let $p$ be a strictly positive element in $B$ and let $a \in X$ such
  that $a^*qa \le p$.  Given $\epsilon > 0$ there exists $b \in X$
  such that $qb = qa$ and $b^*b \le p + \epsilon 1_B$.
\end{proposition}

\begin{proof}First assume $p =1$.  Let $I = \{x \in X : qx = 0\}$ and let
  $\delta>0$ such that $2\delta + \delta^2 < \epsilon$.  Then by the
  previous lemma there exists an $h \in I$ such that $\| a+h \| \le \|
  qa \| + \delta$.  Let $b = a+h$ and note that $qb = qa$.  Also,
  since $a^*qa \le 1$, we have $\|qa\| \le 1$.  Then
  \begin{eqnarray*}b^*b & \le & \|b^*b\|1_B = \|b\|^21_B \\
    & \le & (\|qa\| + \delta)^21_B \le (1 + \delta)^21_B \\
    & = & (1 + 2\delta +\delta^2)1_B \\
    & \le   & (1 + \epsilon)1_B = p + \epsilon 1_B.
  \end{eqnarray*}
  In the case that $p$ is not necessarily $1$, note that $a^*qa \le p$ is
  equivalent to $$p^{-1/2}a^*qap^{-1/2} \le 1.$$  Furthermore, note that
  $p^{-1/2}a^*qap^{-1/2} = (a p^{-1/2})^*q(a p^{-1/2})$.  Now suppose
  that $\varphi \in (Xp^{-1/2})^\perp \subset B^*$.  Then $\varphi(\cdot
  p^{-1/2})\in X^\perp$.  Hence, by hypothesis, $\varphi(\cdot
  p^{-1/2})\in (qX)_\perp,$ and thus $\varphi \in  (qXp^{-1/2})^\perp.$
  So by the $p=1$ case, there exists $bp^{-1/2} \in Xp^{-1/2}$ such
  that 
  $$qbp^{-1/2} = qap^{-1/2},$$
  and
  $$p^{-1/2}b^*bp^{-1/2} \le 1 + \epsilon \| p \|^{-1}.$$  
  Pre- and post- multiplying by $p^{1/2}$ yields
  $$b^*b \le p + \epsilon \|p\|^{-1} p \le p + \epsilon.$$
\end{proof}
 
\begin{proposition}\label{noepsilon}Let $X$ be a
  unital subspace of $B$ and suppose $q$ is a projection in
  $B^{**}$ such that $\varphi \in (qX)_{\perp}$ for every 
  $\varphi \in X^{\perp}$.  Let $p$ be a strictly positive contraction
  in $B$.  If $a \in X$ with $a^*qa \le p$, then 
  there exists $b$ in the unit ball of $\{ x \in X : qx = qa
  \}^{\perp\perp}$ such that $b^*b \le p$.  Moreover, $qb = qa$.  
\end{proposition}

\begin{proof}As in the previous proposition, we first show that the
  lemma holds in the case $p=1$.
  Suppose $p=1$. By the previous lemma, for
  each $n>1$ there is a $b_n \in X$ such that $qb_n = qa$ and $b_n^*b_n \le
  1 + \frac{1}{n}$.  By the weak*-compactness of
  $\text{Ball}(\overline{X}^{w^*})$, 
  $(b_n)$ has a weak*-limit point 
  $b$ in $\text{Ball}(\overline{X}^{w^*}).$ 
  Thus $b^*b \le 1$.  Let $(b_{n_t})$ be a subnet of $(b_n)$
  converging to $b$.  Then by weak*-continuity we must also have that
  $qb = qa$.  For general $p$, as before we note that $a^*qa \le p$ is
  equivalent to $$p^{-1/2}a^*qap^{-1/2} \le 1,$$ and that
  $p^{-1/2}a^*qap^{-1/2} = (a p^{-1/2})^*q(a p^{-1/2})$.  Now 
  let $\varphi \in (Xp^{-1/2})^\perp$.  Then $\varphi(\cdot
  p^{-1/2})\in X^\perp$.  Thus, by hypothesis, $\varphi(\cdot
  p^{-1/2})\in (qX)_\perp,$ and thus $\varphi \in  (qXp^{-1/2})^\perp.$
  So by the $p=1$ case, there exists $bp^{-1/2} \in
  \overline{Xp^{-1/2}}^{w^*} = \overline{X}^{w^*}p^{-1/2}$ such
  that 
  $$qbp^{-1/2} = qap^{-1/2},$$
  and
  $$p^{-1/2}b^*bp^{-1/2} \le 1.$$  We pre- and
  post- multiply by $p^{1/2}$ to get $b^*b \le p.$
\end{proof}

{\bf Remarks} \ 1) The preceding two lemmas have matricial
variants.  For 
instance, the conclusion to Lemma \ref{epsilon} can be generalized to
read `for every strictly positive contraction $p \in M_n(B)$ and $a
\in M_n(X)$ with $a^*(I_n \otimes q)a \le  p$, there exists $b \in
M_n(X)$ such that $(I_n \otimes q)a = (I_n \otimes q)b$ and $b^*b \le
p + \epsilon I_n$.'  Here $I_n$ denotes the identity matrix in $M_n$.

2) \ If $X$ is a
  reflexive unital subspace of $B$ and $q$ is such that $\varphi \in
  (qX)_{\perp}$ for every
  $\varphi \in X^{\perp}$, then for every strictly positive
  contraction $p
  \in B$ with $q \le p$, there exists $a \in \text{Ball}(X)$ such that $qa =
  q$ and $a^*a \le p$. 

\bigskip

Variants of Propositions \ref{epsilon} and \ref{noepsilon} in the
commutative case are related 
   to the subject of `peak interpolation' from the theory of function
   algebras (see e.g. \cite{Gam}).  For example, suppose $q$ above is such that
   $\varphi \in A^{\perp}$ implies that $\varphi(q \cdot)=0$, where we
   view $\varphi(q \cdot)$ as an element of $(qB)^*$.   Now suppose
   that $\psi_0 \in (qA)^{\perp}$, viewing $(qA)^{\perp}$ as a
   subspace of $(qB)^*$.  Now define $\psi \in B^*$ by $\psi (b) =
   \psi_0 (qb)$ for all $b \in B$.  Then $\psi \in A^{\perp}$ and, so
   $\psi (q \cdot ) = 0$ as a functional on $qB$, by hypothesis.  Hence,
   we also have $\psi_0(qb) = 0$ for all $b \in A$.  Thus
   $(qA)^{\perp} = \{0\}$ and so $qA$ is norm dense in $qB$.  However, by
   Lemma \ref{gamlemma1}, $qA$ is norm 
   closed, and so $qA = qB$.  Now let $\epsilon > 0$ and let $p$ be a
   strictly positive element of $B$.  Given $a \in B$ with $a^*qa \le
   p$, by Lemma \ref{epsilon}, there exists $b \in A$ such that
   $qb = qa$ and $b^*b \le p + \epsilon$.

\bigskip

We close this section with several lemmas which are required in the
remainder of the paper.  The first one
describes the weak*-limits of powers of certain types of contractions.
 The last two are useful tools for
generating certain closed projections associated with a contraction.

\begin{lemma}\label{weakconverge}
 Let $B$ be a $C^*$-algebra and let $a$ be a contraction in $B^{**}$.
 Let $q$ be a projection in $B^{**}$ such that 
 \begin{enumerate}
   \item $aq = q$, and 
   \item $\varphi(a^*a) < 1$ for all $\varphi
   \in S(B)$ such that $\varphi (q) = 0$.
 \end{enumerate}
 Then $(a^n)$ and $((a^*a)^n)$ converge weak* to $q$ as $n \rightarrow
 \infty$.
\end{lemma}

\begin{proof}
 We have that $B$ is contained non-degenerately in $B(H)$, where
 $H$ is the Hilbert space associated with the universal representation
 of $B$.  We may also view $B^{**}$ as a von Neumann algebra in $B(H)$.
 Let $K$ be the range of $q$ so that $B(H) = B(K \oplus
 K^{\perp})$.  With respect to this decomposition we may write
 $$a = \left[\begin{array}{cc}I_K & 0 \\ 0 & x \end{array} \right],$$
 where $I_K$ is the identity operator on $K$ and $x \in B(K^{\perp})$.  Let
 $\xi \in K^{\perp}$ be a unit vector.  Let $\varphi$ be the vector
 state corresponding to $0 \oplus \xi$.  Then $\varphi (q) = 0$, so that
 $\varphi(a^*a) < 1$.  Thus $\langle x \xi , x \xi \rangle < 1$ for any
 unit vector $\xi$ in $K^{\perp}$.  If $x$ had a reducing subspace on
 which $x$ acted unitarily, then there would be a unit vector $\eta
 \in K^{\perp}$ such that $\langle x \eta , x \eta \rangle = 1$, which
 is a contradiction.  Thus $x$ must be completely non-unitary.  This
 implies that $x^n \rightarrow 0$ in the weak 
 operator topology as $n \rightarrow \infty$.  Now let $\eta_1$ and
 $\eta_2$ be vectors in $K$ and let $\xi_1$ and $\xi_2$ be vectors in
 $K^{\perp}$.  Then,
 \begin{eqnarray*}
   \langle a^n (\eta_1 \oplus \xi_1) , (\eta_2 \oplus \xi_2) \rangle & = & 
   \langle (\eta_1 \oplus x^n \xi_1),(\eta_2 \oplus \xi_1) \rangle \\
   & = & \langle \eta_1 , \eta_2 \rangle + \langle x^n \xi_1 , \xi_2
   \rangle \rightarrow \langle \eta_1 , \eta_2 \rangle.
 \end{eqnarray*}
 Thus $(a^n)$ converges to $q$.  In order to show that $((a^*a)^n)$
 also converges to $q$, it will suffice to show that $x^*x$ is also
 c.n.u.  For in this case, the same argument as above will also work.
 Suppose $x^*x$ has a reducing subspace $V$ on which it acts unitarily
 and let $\xi \in V$ be a unit vector.  Then $\| x \xi \|^2 = \langle x^*x \xi , \xi
 \rangle < 1$, which contradicts $x^*x$ acting unitarily on $V$.
\end{proof}

\begin{lemma}\label{almost_between}Let $X$ be a unital subspace of $B$.  Let
  $a\in \text{Ball}(X)$ and let $q$ be a closed projection in $B^{**}$
  with $aq = q$.  Define $b = \frac{1}{2}(a+1)$.  Then there exists a
  closed projection $r \in B^{**}$ such that $q \le r$  and satisfying
\begin{enumerate}
\item $br=r,$ and
\item $\varphi(b^*b)<1$ for all $\varphi \in S(B)$ such that
  $\varphi(r) = 0.$
\end{enumerate}
\end{lemma}

\begin{proof}With $b = \frac{1}{2}(a + 1)$ and $aq=a$, it is clear
  that $bq = q$.  This
  implies that $(1-b)(1-q) = 1-b$.  Hence, $1-b \in B^{**}(1-q) \cap B$
  and this contains the intersection, $J$, of all left ideals in $B$
  containing $1-b$.  Let $p \in B^{**}$ be the support projection for the
  left ideal $J$.  Then $(1-b)(1-p) = 0$, so that $b(1-p) = 1-p$.
  Hence $1-p$ satisfies condition (1).  Now
  suppose that $\varphi$ is a state of $B$ such that
  $\varphi(1-p) = 0$.  Then surely $\varphi(b^*b) \le 1$, but suppose
  that   $\varphi(b^*b) = 1$.  Then $\varphi(a^*a)+ 2\text{Re} \ \varphi(a) + 1
  = 4$, which 
  forces $\varphi(a^*a) = \varphi(a) = 1$, and hence, $\varphi(b) =
  1$.  Now let $L$ denote the left 
  kernel associated with $\varphi$.  Then $L$ is a left ideal and we
  claim that $1-b \in L$.  To see this, note that
   $$\varphi((1-b)^*(1-b)) = \varphi(b^*b) - 2\text{Re} \;\varphi(b) +1 =
  0.$$  Hence, $1-b \in L$ and consequently, $J \subset L$.
  If $p_L$ denotes the support projection of $L$, then by the
  definition of $J$, we must
  have $p \le p_L$ and so $1-p_L \le 1-p$.  Applying $\varphi$ to this
  last inequality yields $1 \le 0$, an obvious contradiction.  Thus
  we conclude that $ \varphi(b^*b)  < 1$ and  so $1-p$ satisfies
  condition (2).  Furthermore, since $J \subset B^{**}(1-q)\cap
  B$, it follows that $q \le 1-p$.  Now let $r = 1-p$.
\end{proof}

We also need the following variant of Lemma \ref{almost_between}:

\begin{lemma}\label{weak*_between}Let $a$ be a contraction in $B^{**}$
  and $q$ a closed 
  projection in $B^{**}$ such that $aq = q$.  Let $b = \frac{1}{2}(a+1)$.
  Then there exists a projection $r \in B^{**}$ such that $r \ge q$
  and satisfying
  \begin{enumerate}
   \item $br = r$, and
   \item $\varphi(b^*b) < 1$ for all $\varphi \in S(B)$ such that
     $\varphi(r)=0$.
  \end{enumerate}
  Moreover, $b^k \rightarrow r$ weak* as $k \rightarrow \infty$.
\end{lemma}

\begin{proof}This is essentially the same proof as above.  Let $b =
  \frac{1}{2}(a+1)$ and let $J$ be the intersection of all
  weak*-closed left ideals in $B^{**}$ containing $1-b$.  This will be
  a weak*-closed left ideal.  Let $p$ be the support
  projection of $J$.  Then, as above, $q \le 1-p$ and $b(1-p) = 1-p$.
  Now let $\varphi \in S(B)$  be such that $\varphi(1-p) = 0$.  Letting
  $L$ be the left-kernel of $\varphi$ in $B^{**}$, we get a
  weak*-closed left ideal containing $1-b$.  As in the proof of Lemma
  \ref{almost_between}, we let $p_L$ be the support projection of $L$,
  so that $1-p_L \le 1-p$.  Applying $\varphi$ to this inequality
  yields the necessary contradiction.  Taking $r = 1-p$ proves
  the first part.  The second part follows from Lemma \ref{weakconverge}.
\end{proof}

\section{Right ideals with left contractive approximate identity}\label{right}

We now have everything needed to prove the main theorem.  The
following theorem gives the difficult direction. 

\begin{theorem}\label{qinclosure}Let $A$ be a unital subalgebra of a
  unital $C^*$-algebra 
  $B$.  Let $q\in B^{**}$ be a closed projection such that $q \in
  A^{\perp\perp}$.  Then $1-q$ is in the weak*-closure of the right
  ideal $J = \{a \in A :(1-q)a = a \}$.
\end{theorem}

\begin{proof}
  First recall that $A^{\perp} \subset (qA)_{\perp}$ is equivalent to
  $q \in A^{\perp\perp}$ by Lemma \ref{weak}.  Thus $q$ satisfies the
  hypotheses of Proposition \ref{epsilon} and \ref{noepsilon}.  
  Let $u$ be an open projection dominating $q$.  Then, by
  the noncommutative Urysohn's lemma there exists $p \in B$, $0 \le
  p \le 1$, such that $pq = q$ and $p(1-u) = 0$.  For each integer $n
  \ge 0$, let $$p_n = \frac{n}{n+1}p + \frac{1}{n+1}.$$  Then $p_n$ is
  strictly positive and $p_nq = q$, so that $q \le p_n$.  Each $p_n$
  also has the property that $p_n(1-u) = \frac{1}{n+1}(1-u)$.  By Proposition
  \ref{noepsilon}, 
  for each $n$ there exists $a_n \in \text{Ball}(A^{**})$ such that
  $qa_n = q$ and 
  $a_n^*a_n \le p_n$.  The net $(a_n)$ is contained in the
  unit ball of $A^{**}$, which is weak*-compact.  Let $(a_{n_t})$ be a
  subnet converging to an element $a$ in the unit ball of $A^{**}$.
  Since $(1-u)a_{n_t}^*a_{n_t}(1-u) \le \frac{1}{n_t+1}$, the 
  net $(1-u)a_{n_t}^*a_{n_t}(1-u)$ converges to zero in norm, and hence,
  by the $C^*$-identity,
  $a_{n_t}(1-u)$ also converges to zero in norm. It follows that $a(1-u) =
  0$ and $au = a$. 
  Similarly, $qa = q$.  Let $b =
  \frac{1}{2}(a+1)$.  From Lemma \ref{weak*_between} we know that there is a
  projection $r \in B^{**}$ with $r \ge q$ such that $b^k \rightarrow r$ weak*.
  We now show that $r \le u$.  To do this, we first observe that
  \begin{eqnarray*} bu &  = &  \left(\frac{1}{2}a +
  \frac{1}{2}\right)u \\
    & = & \frac{1}{2}a + \frac{1}{2}u \\
    & = & \frac{1}{2}a + \frac{1}{2}u +  \frac{1}{2}(1-u) -
    \frac{1}{2}(1-u) \\
    & = & b-\frac{1}{2}(1-u),
  \end{eqnarray*}
and so $$b^k u =  b^k - \frac{1}{2}b^{k-1}(1-u).$$

Thus in the weak*-limit, as $k \rightarrow \infty$, we get
 $$ru  =  r - \frac{1}{2}r(1-u).$$
Hence, $$2ru = 2r -r(1-u) = r + ru,$$ and therefore $ru = r$.

For each $a_n$, Lemma \ref{weak*_between} gives rise to
projections $q_n$ in $B^{**}$ with $q_n \ge q$ such that $b_n \equiv \frac{1}{2}(a_n+1)$ has
the following properties:\begin{enumerate}
\item $q_nb_n = q_n$,
\item $\varphi(b_n^*b_n) < 1$ for all $\varphi \in S(B)$ such that
  $\varphi(q_n)=0$, and
\item $b_n^k \rightarrow q_n$ weak* as $k \rightarrow
  \infty$.\end{enumerate} 
Item (1) implies that $1-b_n$ is in $J$ and item (3) implies that each
  $q_n$ lies in $A^{**}$.

Now let $Q = \bigwedge q_n$, so that $Q \le q_n$ for all $n$.  Since $q_n
\ge q$ we also have $Q \ge q$.  The containment $B^{**}(1-q_n)
\subset B^{**}(1-Q)$ follows from $Q \le q_n$.  However, $1-b_n \in
B^{**}(1-q_n)$ for each $n$.  Hence $1-b_n \in B^{**}(1-Q)$ for all
$n$.  However, $1-b_{n_t} \rightarrow 1-b$, and so $1-b$ is in
$B^{**}(1-Q)$.  By the construction of $r$ (recall $1-r$ is the support
projection for the weak*-closed left ideal generated by $1-b$) this
implies that $B^{**}(1-r) 
\subset B^{**}(1-Q)$.  Therefore $q \le Q \le r \le u$ and hence,
$$q \le \bigwedge_n q_n \le u.$$  Let $Q_u =
\bigwedge_n q_n$, which is in $A^{**}$. As $u$ varies over
all open projections dominating $q$, we get
$$q \le \bigwedge_{u \ge q}Q_u \le \bigwedge_{u \ge q}u = q,$$
So, $$q = \bigwedge_{u \ge q}Q_u =
  \bigwedge_{u \ge q} \bigwedge_n q_n.$$  Thus,
  \begin{eqnarray*}1-q & = & 1-  \bigwedge_{u \ge q} \bigwedge_n
    q_n \\
    & = & \bigvee_{u \ge q}\left( 1-\bigwedge_n q_n \right) \\
    & = & \bigvee_{u \ge q}\bigvee_n \left( 1- q_n \right).
  \end{eqnarray*}


  If $u$ is fixed, then for each $q_n$ associated with $u$, $b_nq =
  q$, where $b_n \in A^{**}$, as above.  From Proposition \ref{noepsilon},
  each $a_n$ is a weak*-limit of elements $y$ in $A$ satisfying
  $qy=q$.  So each $b_n$ is the
  weak*-limit of a net, $(c_t)$ say, in $A$ such that $qc_t = q$.  Thus
  $(1-q)(1-c_t) = 1-c_t$ and therefore the net $(1-c_t)$ is contained in
  $J$.  Hence its weak*-limit, $1-b_n$, is in $J^{\perp\perp}$.  Also,
  for any integer $k>0$, $1-b_n^{k}$ is in $J^{\perp\perp}$.  Hence
  $1-q_n$ = $w^*$-$\lim_k 1-b_n^{k}$ is in $J^{\perp\perp}$.  Combining
  this last fact with the last displayed equation we see that $1-q$ is
  in $J^{\perp\perp}$.
\end{proof}

As a consequence, we have our main theorem characterizing right ideals
with left contractive approximate identity.

\begin{theorem}Let $A$ be a unital subalgebra of a unital
  $C^*$-algebra $B$.  A subspace $J$ of $A$ is a right ideal
  with left contractive approximate identity if and only if $J =
  (1-q)A^{**} \cap A$ for a closed projection $q \in A^{\perp \perp}$.
\end{theorem}

\begin{proof}The forward implication is the easy direction and is
  essentially in \cite{BEZ}.  Suppose
  $J$ is a right ideal with left contractive 
  approximate identity $(e_t)$.  Then $J^{\perp \perp}$ has a left
  identity $p$ such that $e_t \rightarrow p$ weak* (see e.g. 2.5.8 in
  \cite{BL}).  Since $p$ is a contractive idempotent, it is an
  orthogonal projection.  Since $pe_t = e_t$, $p$ is an open projection
  by the discussion on open projections in Section \ref{prelim}. 
  So $q = 1-p$ is closed.  Also, $J \subset (1-q)A^{**} \cap A$.
  However, if $a \in A$ such that $(1-q)a = a$, then $e_ta \in J$ and
  so $a = (1-q)a \in J^{\perp \perp} \cap A = J$.  Thus 
  $(1-q)A^{**} \cap A \subset J$.  Now let $\varphi \in A^{\perp}$.  Then $0 =
  \varphi(1-e_t) \rightarrow \varphi(q)$, and so $\varphi (q) = 0$.
  Hence $q \in A^{\perp \perp}.$

  On the other hand, suppose $J$ is a subspace of $A$ such that $J =
  (1-q)A^{**} \cap A$ for a closed projection $q \in A^{\perp
  \perp}$.  The subspace $J$ is
  a right ideal of $A$ and $J = \{a \in
  A : qa = 0 \}$.  By Theorem \ref{qinclosure}, $J^{\perp \perp}$
  contains $1-q$, which is a left identity for $J^{\perp \perp}$.
  Thus $J$ possesses a left contractive approximate identity (see
  e.g. 2.5.8 in \cite{BL}).
\end{proof}

\section{peak projections}\label{peak}

From the last section we see can see the role that closed projections
in $A^{\perp\perp}$ play in determining the right ideal structure of
an operator algebra $A$.  In this section we study the
closed projections in $A^{\perp\perp}$ from the view
point of `peak phenomena' in $A$.  Indeed, this idea was already subtly
playing a role in the proof of Theorem \ref{qinclosure}.
The following theorem will be the basis for our definition of a
noncommutative peak set, or \emph{peak projection}. However, first
note that if $a$ is a contraction and $q$ is a projection such that
$qa = q$, then $a$ and $q$ necessarily commute.

\begin{theorem}\label{equivalent}Let $a$ be a contraction in $B$
  and let $q$ be a 
  closed projection in $B^{**}$ such that $aq=q$.  Then the following
  are equivalent:
  \begin{enumerate}
   \item $\varphi(a^*a) < 1$ for all $\varphi \in S(B)$ such that
     $\varphi(q) = 0$, 
   \item $\varphi (a^*a(1-q)) < 1$ for all $\varphi \in S(B)$,
   \item $\varphi (a^*a(1-q)) < \varphi (1-q)$ for all $\varphi \in
     S(B)$ such that $\varphi(1-q) \not= 0$,
   \item $\varphi(a^*a) < 1$ for every pure state $\varphi$ of $B$ such
     that $\varphi(q) = 0$,
   \item $\| pa \| < 1$ for any closed projection $p \le 1-q$,
   \item $\| ap \| < 1$ for any closed projection $p \le 1-q$, and
   \item $\| ap \| < 1$ for any minimal projection $p \le 1-q$. 
  \end{enumerate}
\end {theorem}

\begin{proof} (3) $\Rightarrow$ (2) \  Assume (3) holds.  If $\varphi \in
  S(B)$ is such that $\varphi$ doesn't 
  vanish on $1-q$, we have $\varphi (a^*a(1-q)) <
  \varphi (1-q) = 1 - \varphi(q) \le 1$.  In the case that $\varphi$
  vanishes on $1-q$, (2) follows by the Cauchy-Schwarz inequality for
  positive linear functionals.  
  
  (2) $\Rightarrow$ (3) \ Suppose (2) holds.  Let $\varphi$ be a state
  which does not vanish on $1-q$ and 
  define $$\psi ( \cdot ) = \frac{\varphi (\cdot (1-q))}{\varphi
    (1-q)}.$$  
  Because $\psi$ is contractive and unital, it is a state on $B$.
  Applying $\psi$ to both sides of
  $a^*a(1-q) < 1$ we get
  $$\frac{\varphi (a^*a(1-q)(1-q))}{\varphi(1-q)} < 1,$$
  which implies that $\varphi (a^*a(1-q)) < \varphi (1-q).$ So (3) holds.

  (1) $\Rightarrow$ (4) \  This is immediate.

  (2) $\Rightarrow$ (1) \  Assume (2) holds and let $\varphi \in S(B)$ be
  a state such that 
  $\varphi(q) = 0$.  Then $\varphi(a^*a) = \varphi(a^*a) - \varphi(q) =
  \varphi(a^*a(1-q)) < 1$.  So (1) holds.

  (4) $\Rightarrow$  (2) \ Assume (4) and let $\varphi \in S(B)$ and
  suppose that 
  $\varphi(a^*a(1-q)) = 1$.  It follows that $\varphi(a^*a) = 1$
  and $\varphi(q) = 0$.  Consequently, $\varphi(1- a^*a) = 0$.
  Now let $J$ be the left ideal $B^{**}(1-q) \cap B$ of $B$, so that
  $J \cap J^* = (1-q)B^{**}(1-q) \cap B$ is a hereditary subalgebra of
  $B$.  Note that  $1-a^*a = (1-a^*a)(1-q)$ is an element of $J \cap
  J^*$.  Let $(e_t) \subset J \cap J^*$ be an increasing contractive
  approximate identity for $J \cap J^*$.  Then $(e_t)$ is an
  increasing left contractive approximate identity for $J$, However,
  since $1-q$ is the support 
  projection for $J$, then necessarily $e_t \rightarrow 1-q$ weak*, and so
  $\varphi(e_t) \rightarrow \varphi(1-q) = 1$.  However,

  $$\|\varphi|_{J \cap J^*} \| = \lim_t \varphi|_{J \cap
  J^*}(e_t)$$ and so  $\|\varphi|_{J \cap J^*} \| = 1$, making
  $\varphi|_{J \cap J^*}$ a state of $J \cap J^*$ such that
  $\varphi|_{J \cap J^*}(1-a^*a) = 0$.  By Lemma \ref{vanishing}
  there exists a pure state $\psi$ of  $J \cap J^*$ which annihilates
  $1-a^*a$.  We can then extend $\psi$ to a pure state $\tilde{\psi}$
  of $B$.  We then have
  $$0 =  \tilde{\psi}(1-a^*a) = \tilde{\psi}(1) - \tilde{\psi}(a^*a) =
  1 - \tilde{\psi}(a^*a).$$  Thus $\tilde{\psi}(a^*a) = 1$.  We also
  have
  \begin{eqnarray*}1 &  = &  \tilde{\psi}(1) = \tilde{\psi}(q) +
    \tilde{\psi}(1-q)\\
    & = & \tilde{\psi}(q) + \lim_t \tilde{\psi}(e_t) \\
    & = & \tilde{\psi}(q) + \lim_t \psi(e_t) \\
    & = & \tilde{\psi}(q) + \|\psi\| = \tilde{\psi}(q) + 1,
  \end{eqnarray*} which forces $\tilde{\psi}(q) = 0$.  Thus we have
  found a pure state on $B$ which annihilates $q$, but takes the value
  1 at $a^*a$. This contradicts our assumption of (4).  Thus, the
  supposition that $\varphi(a^*a(1-q)) = 1$ must be rejected, and
  therefore (2) holds. 

  (6) $\Rightarrow$ (7) \  This is
    immediate from the fact that any minimal projection is automatically
    closed. 
  
  (7) $\Rightarrow$ (4) \  Assume (7) and let $\varphi$ be a pure state
    of $B$ which 
    annihilates $q$.  Let $\varphi$ be such a pure state.  Let $L$ be
    the left-kernel of $\varphi$. Then $L =B^{**}(1-p)\cap B$ for
    some minimal closed projection $p \in B^{**}$, and the
    weak*-closure of $L$ in $B^{**}$ will be $B^{**}(1-p)$ (3.13.6 in
    \cite{Ped}). Now viewing $\varphi$ as a normal state on $B^{**}$,
    let $L'$ be the left-kernel of $\varphi$ with respect to $B^{**}$.
    This will be a  weak*-closed left ideal and $L' = B^{**}r$ for
    some projection $r \in B^{**}$  with $q \in L'$.  The 
    containment $B^{**}(1-p) \cap B \subset B^{**}r$ is obvious.
    Passing to the weak*-closure we get $B^{**}(1-p) \subset B^{**}r$.
    Hence, $1-p \le r$, or rather $1-r \le p$, which by minimality of
    $p$ forces $p = 1-r$.  We conclude that $\overline{L}^{w^*} = L'$
    and therefore $q \in B^{**}(1-p)$ .  Hence $q \le 1-p$, or
    equivalently, $p \le 1-q$.
    Thus by condition (7), $\|ap\| < 1$.    Now decompose $a^*a$ as
    $$a^*a = pa^*ap + pa^*a(1-p) + (1-p)a^*ap + (1-p)a^*a(1-p)$$ and
    apply $\varphi$ to get
    $$ \varphi (a^*a) = \varphi(pa^*ap) + \varphi(pa^*a(1-p)) +
    \varphi((1-p)a^*ap) + \varphi((1-p)a^*a(1-p)).$$ The last 3 terms
    vanish by the Cauchy-Schwarz inequality and the fact that
    $\varphi (1-p) = 0$.  Thus we have $$\varphi (a^*a) =
    \varphi(pa^*ap) \le \|pa^*ap\| = \|ap\|^2 < 1,$$ establishing (4).

   $(2) \Rightarrow (5)$ \  Assume (2) and let $p \le 1-q$ be closed.
    Since $p$ is closed,  $a^*pa$ is upper semi-continuous on $S(B)$.
    Therefore, $a^*pa$ attains its maximum on $S(B)$ at some state
    $\varphi$, and hence $\|pa\|^2 = \varphi(a^*pa)$.  Since  $p
    \le 1-q$ we have $a^*pa \le a^*(1-q)a$, and so $$\varphi(a^*pa)
    \le \varphi(a^*(1-q)a) =  \varphi(a^*a(1-q)) < 1,$$ by (2) and
    the fact that $a$ and $q$ commute.
    Thus, the norm of $pa$ must be strictly less than 1.

   $(6) \Rightarrow (7)$ \  Apply the implication $(7) \Rightarrow
    (6)$, which has already been established, to  $a^*$.
  \end{proof}

{\bf Remarks} \ First, the condition $\|ap\| < 1$ is equivalent to
$\|pa^*\| < 1$, which means 
that $a^*a$ in conditions (1)-(4) can be replaced with $aa^*$. Second,
similarly to the equivalence of (6) and (7), condition (5) is equivalent to the
statement that $\|pa\| < 1$ for any minimal projection $p \le 1-q$. 

\begin{definition}\label{peakdef}
  If $A$ is a unital subalgebra of a
   unital $C^*$-algebra $B$, a projection $q \in B^{**}$ is called a
  \emph{peak projection for $A$} if there exists a contraction $a
  \in A$ such that $qa = q$ and such that $a$ and $q$ satisfy the
  equivalent conditions of Theorem \ref{equivalent}.
  We refer to $a$ as the \emph{peak associated with $q$}. 
  If $q$ is an intersection of peak projections, we refer to $q$ as  a
  \emph{$p$-projection}.
\end{definition}

If $B = C(K)$, the continuous complex-valued functions on a
compact space $K$, then the peak projections for $A$ are exactly the
characteristic functions of peak sets.  It
is also easy to see that if $E$ is a peak set and $a$ is its associated
peak, then $a^n$ converges point-wise to $\chi_E$.  From Lemma
\ref{weakconverge}, the same is true
for peak projections.

As a consequence, we have the following proposition.

\begin{proposition}Any projection satisfying condition (1) in Theorem
  \ref{equivalent} for some contraction in a unital subspace $A$
  of a unital $C^*$-algebra $B$ is a closed projection.
\end{proposition}

\begin{proof}Let $A$ be a unital subspace of $B$ and let $q$ be a
  projection in $B^{**}$ satisfying condition (1) of Theorem
  \ref{equivalent} for some contraction $a \in A$, then
  $(a^*a)^n$ is a decreasing net of self-adjoint elements in $B$ with
  limit $q$.  Hence $q$ is closed.  Any intersection of closed
  projections is again closed.  Thus the result also holds for
  $p$-projections.
\end{proof}

For a compact Hausdorff space $K$, any closed set in $K$ will be a
$p$-set for $C(K)$ by Urysohn's lemma.  The same thing holds in the
noncommutative case as well. 

\begin{proposition}\label{cstar}For a $C^*$-algebra $B$, any closed projection
  in $B^{**}$ is a $p$-projection for $B$. \end{proposition}

\begin{proof}Let $q \in B^{**}$ be a closed projection.  Then for
  any open projection $u \ge q$, there exists $a_u \in B$ with $0 \le
  a_u \le 1$ such that $a_uq = q$ and $a_u(1-u) = 0$.  Now let $q_u$
  be the weak*-limit of $a_u^n$ as $n \rightarrow \infty$.
  Since multiplication is separately weak*-continuous, it follows that
  $a_uq_u = q_ua_u = q_u$.  From this property and again the separate
  weak*-continuity of multiplication, it follows that $q_u$ is a
  projection.  Since $a_u(1-u) = 0$, we also have $q_u \le u$.  We now
  claim that $q_u$ is a peak projection for $B$ with 
  peak $a_u$.  We only need to check that $\varphi(a_u^2) < 1$ for
  any pure state $\varphi \in B^*$ such that $\varphi(q_u) = 0$.
  However, since   $a_u^n \rightarrow q_u$ weak*, it follows that
  $\varphi(a_u^n) \rightarrow 0$.  Suppose that $\varphi(a_u^2) = 1$.
  Then representing $B$ concretely on a Hilbert space so that
  $\varphi$ is a vector state $\varphi(\cdot) =
  \langle\pi(\cdot)\xi,\xi\rangle$, we see that $\langle \pi(a_u^2)
  \xi, \xi \rangle = 1$.  So by the converse to the Cauchy-Schwarz
  inequality, we must have that $\pi(a_u^2)\xi = \xi$, which
  contradicts that $\varphi(a^n) \rightarrow 0$. So $q_u$ is a peak
  projection. Note also that the equation $a_u^nq = q$ implies that $q
  \le q_u$.

Now we take the intersection $\bigwedge q_u$ of all such $q_u$ as $u$
varies over all open projections dominating $q$.  We now
show that $q = \bigwedge q_u$.  To see this, note that since $q \le
q_u$, we have that $q \le\bigwedge q_u$. By Proposition \ref{regular},
$$q \le \bigwedge q_u 
\le \bigwedge u = q.$$  Thus, $q = \bigwedge q_u.$ \end{proof}

This next proposition describes a peak projection in
terms of a support projection associated with its peak.

\begin{proposition}
A projection $q \in B^{**}$ is a peak projection for a unital subspace
$A$ of $B$ if and only if there exists a contraction $a \in A$ such
that $1-q$ is the right support projection for $1-a$.
\end{proposition}

\begin{proof}We assume that $B$ and $B^{**}$ are acting on the
  universal Hilbert space $H_u$ for $B$.  First suppose that $q$ is a
  peak projection for $A$ with peak $a \in A$.  Let $\xi \in
  H_u$ and let $r$ denote the right support projection of $1-a$. Since
  $a^n \rightarrow q$ weak*,
  $aq = q$ implies that $\xi \in \text{Ran}\;q$ if and only
  if $a\xi = \xi$.  This in turn is equivalent to saying $\xi \in
  \text{Ran}\;q$ if and only if $\xi \in 
  \text{Ker} \; (1-a)$.  Since $r$ is the projection onto $[\text{Ker} \;
  (1-a)]^\perp$, this last statement is equivalent to $1-q = r$. 

  Now suppose that $a \in A$ is a contraction and $q$ is a closed
  projection such that the range projection $r$, of $1-a$, is equal to
  $1-q$.  Then $1-a = (1-a)r = (1-a)(1-q)$ implies $aq=q$. Now
  define $b = (a+1)/2$, so that $bq=q$.  Let $\varphi$ be a state on $B$ such
  that $\varphi(q) = 0$.  We may assume that $\varphi ( \cdot ) =
  \langle \cdot \eta , \eta \rangle$ for a unit vector $\eta \in H_u$.
  Then $q\eta = 0$, and so $r\eta = \eta$.  Now suppose that $\varphi
  (b^*b) = 1$.  Then $$1 = \varphi(b^*b) = \frac{1}{4}(\varphi(a^*a) +
  2\text{Re}\;\varphi(a) + 1).$$  Hence, $1 = \varphi (a^*a) = \langle a \eta
  , \eta \rangle$, and so by the converse to the Cauchy-Schwarz
  inequality, $a \eta = \eta$.  Thus $\xi \in \text{Ker} \; (1-b)$, and
  so $r\eta = 0$, which is a contradiction.  Thus $\varphi (b^*b) < 1$
  and so $q$ is a peak projection with peak $b$.
\end{proof} 

  
In the commutative case, it is easy to see that if $E$ and $F$ are
peak sets for a uniform algebra $A$ with peaks $f$ and $g$,
respectively, in $A$.  then $E \cap F$ will be a peak set with peak
$\frac{1}{2}(f+g)$. For general $C^*$-algebras, we have the following
generalization.

\begin{proposition}Let $q_1$ and $q_2$ be two peak projections with peaks
  $a_1$ and $a_2$, respectively.  Then $q_1 \wedge q_2$ is also a peak
  projection with peak $\frac{1}{2}(a_1 + a_2)$.
\end{proposition}

\begin{proof}
  That  $q_1 \wedge q_2$ and  $\frac{1}{2}(a_1 + a_2)$ satisfy the first
  condition in the definition of peak projection is immediate since $q_1
  \wedge q_2$ is dominated by both $q_1$ and $q_2$.  To show the second
  condition, let $\varphi$ be a pure state of $B$ which annihilates
  $q_1 \wedge q_2$.  Let $b = \frac{1}{2}(a_1 + a_2)$.  We wish to show
  $\varphi(b^*b) < 1$.  This is equivalent to showing
  $$\varphi(a_1^*a_1) +\varphi(a_1^*a_2) + \varphi(a_2^*a_1)
  +\varphi(a_2^*a_2) < 4,$$
  for which it suffices to show that either  $\varphi(a_1^*a_1) < 1$ or
  $\varphi(a_2^*a_2) < 1$.  So suppose $\varphi(a_1^*a_1) =
  \varphi(a_2^*a_2) =1$.  Let $\pi : B^{**} \rightarrow B(H)$ be the
  weak*-continuous cyclic representation associated with $\varphi$ and
  let $\xi \in H$ be 
  the corresponding cyclic vector.  Then $$\langle \pi(a_1^*a_1)\xi ,
  \xi \rangle = \varphi(a_1^*a_1) = 1.$$  Thus, by the converse to the
  Cauchy-Schwarz inequality, we must have $\pi(a_1^*a_1)\xi = \xi$.
  However, since $(a_1^*a_1)^n \rightarrow q_1$ weak*, and so
  $\pi(a_1^*a_1)^n \rightarrow \pi(q_1)$ in the weak operator topology,
  we must have $\varphi(q_1) = 1$.  The same argument shows that
  $\varphi(q_2) = 1$ and so $\xi$ must be in the range of both
  $\pi(q_1)$ and $\pi(q_2)$.  Hence, $\varphi(q_1 \wedge q_2) = \langle
  \pi(q_1) \wedge \pi(q_2)\xi , \xi  \rangle = 1,$ which is a contradiction.
\end{proof}

The next result shows that $p$-projections must be in the
weak* closure of $A$ in $B^{**}$ when $A$ is a subalgebra.

\begin{proposition}\label{forward}Let $A$ be a unital subalgebra of a unital
  $C^*$-algebra $B$, and let $q$ be a $p$-projection for $A$.  Then
  $\varphi(q) = 0$ for all $\varphi$ in $A^{\perp}$.  Consequently,
  $\varphi$ is in $(qA)_{\perp}$ and $q$ is in $A^{\perp \perp}$.
\end{proposition}

\begin{proof}First assume $q$ is a peak projection and let $a \in A$ be
  the peak associated with $q$.  For any  $\varphi$ in $A^{\perp}$ and
  any integer $n>0$, $\varphi (a^n) = 0$.  However, $a^n$ converges
  weak* to $q$.  Thus $\varphi (q) = 0$.  
  
  Now suppose  $q = \bigwedge_iq_i$ and let $\varphi \in A^{\perp}$.
  Let $\epsilon > 0$.  By a result in \cite{Ak1} there exists an open
  projection $p \in B^{**}$ such that $p \ge q$ and $|\varphi|(p-q) <
  \epsilon$, where $|\varphi|$ is obtained from the polar
  decomposition of $\varphi$ (3.6.7 in \cite{Ped}).  By hypothesis
  $\bigwedge_iq_i \le p$.  Hence $1-p \le \bigvee_i(1-q_i)$, and so by
  the compactness property of closed projections there, exist finitely
  many projections 
  $q_1,q_2,\dots,q_n$ in the family $\{q_i\}$ such that $1-p \le
  \bigvee_{i=1}^n (1-q_i)$.  Thus $q \le \bigwedge_{i=1}^n q_i \le p$.  Now
  let $Q = \bigwedge_{i=1}^n q_i$, which is again a peak projection.
  By the last paragraph it 
  follows that $\varphi(Q) = 0$, and so $|\varphi(Q-q)| =
  |\varphi(q)|$ for all $\varphi \in A^{\perp}$.  The functional
  $|\varphi|$ has the property that  
  $|\varphi(x)|^2 \le \|\varphi\|  |\varphi|(x^*x)$ for all $x \in
  B^{**}$. Thus we have
\begin{eqnarray*}|\varphi(Q-q)|^2 & \le & \|\varphi\|  |\varphi|((Q-q)^*(Q-q))\\
  & = & \|\varphi\|  |\varphi|(Q-q))\\
  & \le & \|\varphi\| |\varphi|(p-q)) \\
  & < & \|\varphi\| \epsilon. \end{eqnarray*}
Since $\epsilon$ was arbitrary, this shows that $\phi(q) = 0$.

  If $x \in
  A$, the map $\varphi(\cdot x)$ is also in $A^{\perp}$.  Thus
  $\varphi(qx) = 0$ for every $x$ in $A$. This shows that
  $\varphi$ is in  $(qA)_{\perp}$.  By Lemma \ref{weak}, this implies
  that $q$ is in $A^{\perp\perp}$. 
\end{proof}

It is natural to ask if the
notion of a peak or $p$-projection is dependent on the particular
$C^*$-algebra in which we view $A$ as residing.  That is, if we have 
embeddings of $A$ into two different $C^*$-algebras, can the
peak projections arising from the each embedding be identified in some
way?  By an embedding of $A$ into a $C^*$-algebra, we mean a completely
isometric homomorphism of $A$.  The following proposition shows that
the notion of a peak or $p$-projection is indeed independent of the
particular embedding.

\begin{proposition}Let $A$ be a unital subalgebra of a $C^*$-algebra
  $B$.  Let $\pi:A \rightarrow B_1$ be a unital completely isometric
  homomorphism of $A$ into another $C^*$-algebra $B_1$.  If $q$
  is a $p$-projection for $A$ in $B^{**}$, 
  then $\pi^{**}(q)$ is a $p$-projection for $\pi(A)$ inside
  $B_1^{**}$.\end{proposition}

\begin{proof}We assume that $B_1$ is acting on its universal Hilbert
  space $H_u$, that is $B_1 \subset B(H_u)$.  Assume that $q$ is a peak
  projection with peak $a \in
  A$. Clearly, $\pi(a)\pi^{**}(q) = \pi^{**}(q)$.  Now suppose
  that $\varphi \in S(B_1)$ such that $\varphi(\pi^{**}(q)) = 0$ and
  extend $\varphi$ to a vector state $\tilde{\varphi}$
  on $B(H_u)$.  Viewing $\pi$ as a unital completely positive map into
  $B(H_u)$, by Arveson's extension theorem we may extend $\pi$ to a
  completely positive map $\rho: B \rightarrow B(H_u)$.  This in turn
  can be extended to a weak*-continuous map $\tilde{\rho}$ on
  $B_1^{**}$.  Since
  $\pi^{**}$ is the unique weak*-continuous extension of $\pi$ to
  $A^{**}$, we must have $\tilde{\rho}|_{A^{**}} = \pi^{**}$.  We
  next observe that $\tilde{\varphi}\circ\rho$ is a state on $B$ which
  extends uniquely to a weak*-continuous state on $B^{**}$, which by
  uniqueness must be $\tilde{\varphi} \circ \tilde{\rho}$.  Hence,
  $(\tilde{\varphi} \circ \rho)(q) = \tilde{\varphi}(\tilde{\rho}(q)) =
  \varphi(\pi^{**}(q)) = 0$.  Thus, since $q$ is a peak projection, we
  must have that
  $\tilde{\varphi}(\rho(a^*a)) < 1$.  By the Kadison-Schwarz inequality
  for completely positive maps, we have $$\varphi(\pi(a)^*\pi(a)) =
  \tilde{\varphi}(\rho(a)^*\rho(a)) \le 
  \tilde{\varphi}(\rho (a^*a)) < 1.$$
  Thus $\pi^{**}(q)$ is a peak projection for $\pi(A)$.
  Now if $q$ is just a $p$-projection
  with $q = \wedge q_i$ for peak projections $q_i$, then $q$ is the
  weak* limit of the net of meets for finitely many $q_i$.  Thus, by
  the discussion about meets in Section 2, it follows that
  $\pi^{**}(q) = \wedge 
  \pi^{**}(q_i)$.  Thus $\pi^{**}(q_i)$ is a $p$-projection.
\end{proof}

Minimal projections which are also $p$-projections correspond to
\emph{$p$-points} (an intersection of singleton $p$-sets) in the commutative case.  The closure of the set of
$p$-points for a uniform algebra is the Shilov boundary (see
e.g. \cite{Gam} and \cite{Phelps}).
Let $A$ and $B$ be as before, but assume $A$ generates $B$ as a
$C^*$-algebra.  Then there exists a largest closed
two-sided ideal $J$ of $B$ such that the canonical quotient map $B
\rightarrow B/J$ restricts to a complete isometry on $A$
(\cite{Arv1}).  The ideal $J$ is the 
so-called `Shilov ideal' for $A$.  Let $p$ be the closed
projection in $B^{**}$ corresponding to the Shilov boundary ideal,
then $p$ dominates all minimal projections which are also $p$-projections
for $A$. Moreover, $p$ dominates all minimal projections in $A^{\perp\perp}$ (see \cite{BHN}). 
 
Unfortunately, however, the join of the orthogonal complements of all
such minimal projections does not in general equal the support of the Shilov 
ideal.

\bigskip

One of the interesting aspects of $p$-projections is their
relationship to approximate identities.  For instance, we have the
following proposition (\cite{Hay},\cite{BHN}).

\begin{proposition}\label{lcai}If $A$ is a unital subalgebra of a
  $C^*$-algebra $B$ and $p \in B^{**}$ is the support projection for a
  right ideal of $A$ with a left approximate identity of
  the form $(1-x_t)$ for $\|x_t\| \le 1$,  then $1-p$ is a
  $p$-projection for $A$. 
\end{proposition}

\begin{proof}Let $J$ be a right ideal of $A$ with left
  approximate identity $(e_t)$ with $e_t = 1 - x_t$ with $x_t$ in the
  unit ball of $A$.  Let $p$ be its support projection and
  define $q = 1-p$.  Then $q$ is necessarily closed, $e_t \rightarrow
  1-q$ weak*, $J = \{a \in A : (1-q)a = a\}$, and $(1-e_t)q = q$.
  For each $t$ let $J_t$ be the intersection of all right ideals of $B$
  containing $e_t$.  Then there exists a unique closed projection $q_t$ in
  $B^{**}$ such that $J_t = (1-q_t)B^{**} \cap B$.  By the proof of Lemma
  \ref{almost_between} $q_t$
  is a peak projection with peak $\frac{1}{2}[1+(1-e_t)] = 1-\frac{1}{2}
  e_t$ and such that $q_t \ge q$.
  Now set $r = \wedge q_t$.  We have that $r \ge q$, but suppose $r -
  q \not= 0$.  Then $(r-q)e_t \rightarrow (r-q)(1-q) = r-q$,
  and $r-q \le q_t$.  Thus $(1-q_t)B^{**} \subset (1-(r-q))B^{**}$, so
  that $(e_t) \subset (1-(r-q)) B^{**}$.  Hence $(1-(r-q))e_t = e_t$ for
  each $t$, and so $(r-q)e_t=0$.  However,  $(r-q)e_t \rightarrow r-q$.
  Thus $r-q = 0$ and so $r=q$, making $q$ an intersection of peak
  projections.
\end{proof}

{\bf Remark} \  It can also be shown that if an ideal $J$ of a unital
operator algebra $A$ 
has a left contractive approximate identity, then it has a left
approximate identity of the form $(1-x_t)$, where $x_t \in A$ and
$\lim_t \|x_t\| = 1$ (\cite{BHN}).  Moreover, if we can choose the
$x_t$ in $\text{Ball}(A)$, for
every such ideal, then the $p$-projections are exactly the orthogonal
complements of the support
projections for right ideals with left contractive approximate identity.
\bigskip

It is natural to make the following definition.

\begin{definition}Let $A$ be a unital subalgebra of a unital $C^*$-algebra
  $B$.  A projection $q \in B^{**}$ is said to be an \emph{approximate
    $p$-projection for $A$} if $q$ is closed and $q \in A^{\perp\perp}$.
\end{definition}

The following shows that approximate $p$-projections possess peaking
properties.

\begin{theorem}\label{perp}Let $A$ be a unital subalgebra of a unital
  $C^*$-algebra 
  $B$ and let $q \in B^{**}$ be a closed projection.  Then the
  following are equivalent:
  \begin{enumerate}
    \item $q$ is an approximate $p$-projection,
    \item for every $\epsilon > 0$ and for every open projection $u
      \ge q$, there exists $a \in (1+\epsilon)\text{Ball}(A)$ such
      that $qa=q$ and $\|a(1-u)\| \le \epsilon$, and 
    \item for every $\epsilon > 0$ and strictly positive $p \in B$
      with $p \ge q$, there exists $a \in A$ such that $qa = q$ and
      $a^*a \le p + \epsilon$.
  \end{enumerate}
\end{theorem}

\begin{proof}
(1) $\Rightarrow$ (3) \ This is essentially Lemma \ref{epsilon}.  Let
    $\epsilon > 0$ and let $p \in B$ be a
    strictly positive element of $B$ such that $p \ge q$.  Since $q$
    is in $A^{\perp\perp}$, by Lemma
    \ref{weak}, $q$ satisfies the hypothesis of Lemma \ref{epsilon}.
    Thus there exists $a \in A$ such that $qa=q$ and $a^*a \le p +
    \epsilon$.

(3) $\Rightarrow$ (2) \ Let $\epsilon > 0$ and suppose $u \ge q$ is
open.  Let $\delta = \epsilon^2/2$.  
    As we have seen before, by the noncommutative Urysohn's lemma,
    there exists a strictly positive
    contraction $p \in B$ such that $pq = q$  and $p(1-u) = \delta
    (1-u)$.  By (3) there 
    exists $a \in A$ such that $qa=q$ and $a^*a \le p +
    {\delta}.$  Hence, $(1-u)a^*a(1-u) \le 2\delta(1-u),$ and so
    $\|a(1-u)\| \le \epsilon$.

(2) $\Rightarrow$ (1) \ Let $u \ge q$ be an open projection.  By (2), for each
    natural number $n$ there exists $a_n \in (1+1/n)
    \text{Ball}(A)$ such that $qa_n=q$ and $\|a(1-u)\| \le
    \frac{1}{n}$.  The net $(a_n)$ has a weak* limit point $a \in
    \text{Ball}(A^{\perp\perp})$.  Since $\|(1-u)a_n\| \le
    \frac{1}{n}$ for 
    each $n$, we must also have $\|a(1-u)\| = 0$, and hence
    $au = a$.  Let $b = \frac{1}{2}(a+1)$, which is in
    $A^{\perp\perp}$. By Lemma \ref{weak*_between}, there exists a
    projection $q_u \in A^{\perp\perp}$ such that $q \le q_u$ and $b^k
    \rightarrow q_u$.  As in the proof of Theorem \ref{qinclosure} we also
    that $q \le q_u \le u$ for each $u$.  However, by
    Proposition \ref{regular}, this implies that $q = \wedge q_u$ as
    $u$ varies over all open projections dominating $q$.  However, each $q_u$
    is in $A^{\perp\perp}$.  Thus, so is $q$, by the discussion in
    Section \ref{prelim}. 
\end{proof}

For a uniform algebra $A \subset C(K)$, Glickberg's peak set theorem
says that a closed set $E$ is a $p$-set for $A$ if and only if
$\mu \in A^\perp$ implies $\mu_E \in A^\perp$.  With this and Lemma
\ref{weak} in mind, it is then natural to ask whether or not  
the $p$-projections are precisely the approximate $p$-projections, in
the noncommutative setting.
Certainly any $p$-projection is an approximate $p$-projection.  The
reverse implication holds for unform algebras by the classical
Glicksberg theorem, and it holds when $A = B$ is a $C^*$-algebra by
Proposition \ref{cstar}.  It is also true for operator
algebras which are also reflexive Banach spaces, as the following
simple observation shows.  
In particular, it is true for finite dimensional algebras.

\begin{proposition}Let $A$ be a unital subalgebra of a unital
$C^*$-algebra $B$ such that $A$ is also a reflexive Banach space.  
Let $q \in B^{**}$ be a closed projection. 
The following are equivalent:
\begin{enumerate}
  \item $q$ is a $p$-projection for $A$,
  \item $q$ is an approximate $p$-projection for $A$, and
  \item $q \in A$.
\end{enumerate}
\end{proposition} 

\begin{proof}The implication (1) $\Rightarrow$ (2) follows from Proposition
    \ref{weakconverge} and the fact that $A$ is an algebra.  If (2)
    holds, by reflexivity, $q$ is in $A$, establishing (3).  If $q$ is
    in $A$, then it is trivially a  
    $p$-projection.  So (3) $\Rightarrow$ (1) holds.\end{proof}

Approximate $p$-projections enjoy some of the properties of $p$-projections.
For example, by some observations in Section \ref{prelim}, if $A$ is a
unital subalgebra of a unital $C^*$-algebra 
$B$,  then the meet of a collection of approximate
$p$-projections for $A$ is also an approximate $p$-projection.  It can
also be shown that if the join of a collection of approximate
$p$-projections happens to be closed, then it is also an approximate
$p$-projection (see \cite{Hay}).
By Corollary 5.5 of \cite{BHN}, if $q \in B^{**}$ is a closed projection and 
$X$ is a unital subspace of $B$ such that for every strictly positive
contraction $p \in B$ with $q \le p$ there exists $a \in X$ such
that $aq = q$ and $a^*a \le p$,  then $q$ is a $p$-projection for $X$.

As described above, for many algebras the class of $p$-projections
is the same as the class of approximate $p$-projections.  The
most tantalizing remaining question here is whether or not these two notions are the same for a general unital operator algebra.  Nonetheless, the
correspondence between right ideals with left approximate identity 
and approximate $p$-projections will be key in importing some results
from $C^*$-algebra theory to general operator algebras.  Indeed, we
have begun this in \cite{BHN}.


\begin{thebibliography}{99}
\bibitem{Ak1} C.A. Akemann, The general Stone-Weierstrass
    problem, \emph{J. Functional Analysis} {\bf 4} (1969), 277-294.

\bibitem{Ak2} C.A. Akemann, Left ideal structure of
    $C^*$-algebras, \emph{J. Funct. Anal.} {\bf 6} (1970), 305-317.

\bibitem{Arv1} W.B. Arveson, Subalgebras of $C^{*}$-algebras,
  \emph{Acta Math.} {\bf 123 }(1969), 141-224. 

\bibitem{B1} D.P. Blecher, The standard dual of an operator
    space, \emph{Pacific J. Math.} {\bf 153} (1992), 15-30.

\bibitem{B2} D.P. Blecher, One-sided ideals and approximate
    identities in operator algebras, \emph{J. Australian Math. Soc.}
    {\bf 76} (2004), 425-447.

\bibitem{BEZ}D.P. Blecher, E.G. Effros, and V. Zarikian, One-sided
  $M$-ideals and multipliers in operator spaces, I, \emph{Pacific
    J. Math.} {\bf 206} (2002), 287-319.

\bibitem{BHN} D.P. Blecher, D.M. Hay, and M. Neal, Hereditary
    subalgebras of operator algebras, \emph{J. Operator Theory}, to appear.

\bibitem{BL} D.P. Blecher and C. Le Merdy, \emph{Operator Algebras and
  their Modules}, Oxford Univ. Press, 2004.

\bibitem{BZ} D.P. Blecher and V. Zarikian, The calculus of
    one-sided $M$-ideals and multipliers in operator spaces, to
  appear, Mem. Amer. Math. Soc. {\bf 842} (2006).

\bibitem{Dun} N. Dunford and J. Schwartz, \emph{Linear Operators. I.}
  New York: Interscience. 1958.

\bibitem{ER} E. Effros and Z.J. Ruan, \emph{Operator Spaces}, Oxford
  Univ. Press, 2000.

\bibitem{FN} C. Foia\c{s} and B. Sz.-Nagy, \emph{Harmonic Analysis of
    Operators on Hilbert Space}, North-Holland Publishing Company,
  Amsterdam, 1970.

\bibitem{Gam} T. W. Gamelin, \emph{Uniform Algebras}, Prentice-Hall (1969).

\bibitem{G} I. Glicksberg, Measures orthogonal to algebras and
    sets of antisymmetry, \emph{Trans. Amer. Math. Soc.} {\bf 105}
    (1962), 415-435. 

\bibitem{GK} R. Giles and H. Kummer, A non-commutative
    generalization of topology, \emph{Indiana University Mathematics
    Journal}, Vol. 21, No. 1 (1971).

\bibitem{Hay} D.M. Hay, \emph{Noncommutative topology and operator
    algebras}, Ph.D. thesis, University of Houston (in preparation).

\bibitem{Hirsberg} B. Hirsberg, $M$-ideals in complex function
    spaces and algebras, \emph{Israel J. Math.} {\bf 12} (1972), 133-146. 

\bibitem{Jar}K. Jarosz, A characterization of weak peak sets for
  function algebras, \emph{Bul. Austral. Math. Soc.} {\bf 29} (1984), 129-135. 

\bibitem{Kubrusly}C. S. Kubrusly, \emph{An Introduction to Models and
  Decompositions in Operator Theory}, Birkhauser, Boston, 1997.

\bibitem{Paulsen}V. I. Paulsen, \emph{Completely Bounded Maps and
  Operator Algebras}, Cambridge University Press (2002).

\bibitem{Ped}G. Pedersen, {\em $C^*$-algebras and their Automorphism Groups,}
  Academic Press (1979).

\bibitem{Pisier} G. Pisier, \emph{Introduction to Operator Space
  Theory}, Cambridge University Press (2003).

\bibitem{Phelps}R. R. Phelps, \emph{Lectures on Choquet's Theorem},
  Springer-Verlag, Berlin, 2001.

\end{thebibliography}
\end{document}